# Equivalence of L0 and L1 Minimizations in Sudoku Problem

Linyuan Wang*, Wenkun Zhang*, Bin Yan[†], and Ailong Cai

*Abstract*—Sudoku puzzles can be formulated and solved as a sparse linear system of equations. This problem is a very useful example for the Compressive Sensing (CS) theoretical study. In this study, the equivalence of Sudoku puzzles L0 and L1 minimizations is analyzed. In particular, 17-clue (smallest number of clues) uniquely completable puzzles with sparse optimization algorithms are studied and divided into two types, namely, type-I and –II puzzles. The solution of L1 minimization for the type-I puzzles is unique, and the sparse optimization algorithms can solve all of them exactly. By contrast, the solution of L1 minimization is not unique for the type-II puzzles, and the results of algorithms are incorrect for all these puzzles. Each empty cell for all type-II puzzles is examined. Results show that some cells can change the equivalence of L0 and L1 minimizations. These results may be helpful for the study of equivalence of L0 and L1 norm minimization in CS.

*Index Terms*—L1 minimization, equivalence analysis, uniqueness test, Sudoku.

## I. Introduction

Sudoku is a traditional yet interesting puzzle often used as a typical example for optimization algorithm. Each puzzle is presented on a 9×9 square grid in which some digits have already been filled, and the initial occupied cells are called "clues". The aim of solving this puzzle is to fill the empty cells with digits from 1 to 9 such that each digit appears only once in each row, each column, and each 3×3 box. Fig. 1 shows an example of a typical Sudoku puzzle.

Various computer algorithms attempt to solve Sudoku puzzles. A wide range of deterministic algorithms has been proposed based on backtracking, set covering methods, and brutal force search [1–3]. However, because of NP-completeness, most the well-known algorithms in this category have a complexity that expands exponentially with puzzle size and are therefore difficult to solve in general.

Considerable research has applied optimization tools to design low complexity algorithms. In [4], Sudoku puzzle is expressed as a linear system of equations based on the connections with sparse solution. Among the existing optimization algorithms, Gurobi [5], CVX [6], and YALL1 [7] have achieved several fairly good results when puzzles are converted into a mathematical programming model such as L0 and L1 minimization problems. Nonetheless, $P_0$ problem can hardly be solved in general cases [8]. Following the works of [9–11], we convert L0 norm minimization into L1 norm minimization by relaxing the objective convexly that can be easily solved. This category of algorithms seeks to determine a near-optimal approximate "solution" and commonly attains the sparsest solution [12]. However, a problem arises when the L0 and L1 norm minimizations in the context of Sudoku are equivalent. Many theories have been known as sufficient conditions for checking the equivalence of L0 and L1 norm minimizations, including the restricted isometry property [10], the Kashin–Garnaev–Gluskin inequality [11], and the null space property [13]. Nevertheless, these theories can hardly be verified for Sudoku puzzles.

For Sudoku puzzles, there are still many theoretical problems to be discussed, which are closely related to the CS theory. One problem is when the L0 and L1 norm minimizations are equivalent in the context of Sudoku. The study about this should be considered as a helpful example for the study of equivalence of L0 and L1 norm minimizations for the fixed matrix in CS. Another problem concerned is whether the algorithm of L1 minimization can obtain the correct solution when the solution of problem is not unique. This problem is also arisen in other L1 relaxation problems in CS applications.

For the aforementioned theoretical problems in Sudoku puzzles, we primarily aim to verify the equivalence of their L0 and L1 norm minimizations for Sudoku puzzles and obtain the correct answer with L1 relaxation problem in this study.

Fig. 1. Example of a 17-clue Sudoku puzzle

## II. L1-minimization

### A. L1-minimization relaxation

In [4], 9×9 Sudoku puzzles are formulated as a linear system of equations as

$$A * x = b. \qquad (1)$$

[†] Corresponding author is Bin Yan. Email: ybspace@hotmail.com
National Digital Switching System Engineering & Technological R&D Center, Zhengzhou, 450002, China.
* These authors contributed equally to this work.

where $x \in \{0,1\}^{729}$, and $A$ denotes the matrix with different constraints on $x$, $b = [1,1,\cdots,1]^T$.

Given that the solution to the puzzle is sparse, we solve it by solving the below problem [14][15].

$$P_0: \quad \min_x \|x\|_0 \quad s.t.\ A*x = b. \quad (2)$$

Several methods are used to solve L0 minimization directly [16–18], but they usually provide local minimizer to the original problem only.

For $x \in \{0,1\}^{729}$, $P_0$ is equivalent to the problem

$$\min_x \|x\|_1 \quad s.t.\ A*x = b. \quad (3)$$

The preceding problem can be solved as a convex optimization by relaxing the variables as $x \in [0,1]^{729}$.

$$P_1: \quad \min_x \|x\|_1 \quad s.t.\ A*x = b. \quad (4)$$

Similar to the findings in the sparse representation literature [9–11], $P_1$ solves most Sudoku puzzles and identifies the sparsest $x$ that solves (1).

### B. Results of L1-minimization

We use 49151 Sudoku puzzles given in [19] to test our algorithm. A total of 17-clue uniquely completable Sudoku puzzles are already identified, whereas 16-clue examples are yet to be known [20]. [19] collects as many 17-clue examples as possible and obtains 49151 puzzles. We intend to determine the smallest number of entries in a Sudoku puzzle that has a unique completion. Thus, we test all 17-clue uniquely completable Sudoku puzzles.

All test codes are written and tested in MATLAB v7.12.0 run on a Microsoft Windows 7 PC ×64 with Intel(R) Core(TM) i3-4150 CPU 3.50 GHz and 4.00 GB of memory.

CVX and YALL1 are used to solve $P_1$. The results of these methods are the same. They both solve the same 41722 puzzles correctly and fail for the rest. These results are also the same as in [21].

Because of convex relaxation, the incorrect solution got by L1 minimization algorithm contains more than 81 nonzero entries, many of which are equal to 0.5. To our best knowledge, we can not project it onto the correct vector $x \in \{0,1\}^{729}$.

## III. EQUIVALENCE OF L0 AND L1

### A. Uniqueness tests

Previous results inspired us to examine the equivalence of $P_0$ and $P_1$ for every puzzle, that is, the solution of $P_1$ is either unique or not. The following theorem gives necessary and sufficient conditions for a vector $x^*$ to solve $P_1$ uniquely.

**Theorem 3.1 ([22]).** Let $A \in \mathbb{R}^{m \times n} (m < n)$, $I = \mathrm{supp}(x^*) \subseteq \{1,...,n\}$, and $x^*$ is the unique solution of

$$P_1: \quad \min_x \|x\|_1 \quad s.t.\ A*x = A*x^* \quad (5)$$

only if

$$\ker(A_I) = \{0\} \quad (6)$$

and $\omega \in \mathbb{R}^m$ when

$$A_I^T \omega = \mathrm{sign}(x^*), \quad \|A_{I^c}^T \omega\|_\infty < 1. \quad (7)$$

The first condition can be tested by evaluating whether $A_I$ has a full column rank, and the second condition can be examined by converting it into the following optimization problem:

$$\min_{t,\omega} \quad \|A_I^T \omega - \mathrm{sign}(x^*)_I\|_2$$
$$\text{subject to} \quad -t \leq A_{I^c}^T \omega \leq t. \quad (8)$$

$A_I$ is easily observed to have a full column rank because each Sudoku puzzle is uniquely completable. We let $t = 1-\varepsilon$ $(\varepsilon > 0)$ to ensure that the inequality constraint in (7) is satisfied strictly. This constraint can be relaxed and tightened by changing the numerical value of $\varepsilon$.

To study the solution uniqueness of all Sudoku puzzles, we obtain the correct answers for all of them by solving $P_0$ with Gurobi [5][23]. We then let $\varepsilon = 1e-4$ to guarantee that the solution satisfies $\|A_{I^c}^T \omega\|_\infty < 1$ strictly and minimizes the object $\|A_I^T \omega - \mathrm{sign}(x^*)_I\|_2$.

The histogram of the object function values in optimization problem (8) for 49151 puzzles is shown in Fig. 2.

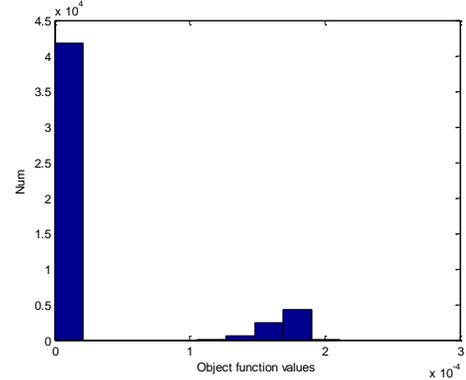

Fig. 2. Histogram of the object function values

The results show that all 41722 puzzles have extremely small object values, which are smaller than 1e-8. Contrarily, the object values of the other puzzles are larger than 1e-4. All the 17-clue puzzles are divided into two different types (i.e., types-I and –II). For the type-I puzzles, the solution of L1 minimization is unique, and the sparse optimization algorithms can solve all of them exactly. For the type-II puzzles, the solution of L1 minimization is not unique, and the algorithms can obtain another solution of $P_1$, which is not the correct answer for the puzzles.

The results above means that researchers could test their improved L1 minimization algorithms for the type-II puzzles. Because the $P_0$ problem is NP-hard, all the improved L1 minimization algorithms may not solve all the type-II puzzles correctly.

There are some interesting theoretical questions in our results. We note that the size of the Sudoku system matrix for all the 17-clue puzzles are 341×729, and 324 rows are exactly same. Different 17 rows come from 17 initial clues. What is the difference of the matrix between the type-I puzzles and type-II? What structure of property is natural? Are there some new conditions like RIP for the Sudoku matrix? We think these open questions are very useful for the study of equivalence in CS.

*B. Further study*

Then, we study the equivalence further. Each empty cell of all the type-II puzzles is examined by testing the uniqueness of the 18-clue puzzle obtained by filling the true number of one empty cell for a 17-clue puzzle. We test the uniqueness for all 475456 (64×7429) puzzles for 18-clue. Experiment results show that some cells can change the equivalence of L0 and L1 minimizations. This finding implies that if we obtain the value of these cells, then the puzzle can be correctly solved. Fig. 3 shows the position of the key cells for the 17-clue puzzle in Fig. 1 with red color. No rule explains the position of these key cells. However, we find most of the key cells have the same value in the CVX solution whose correct value is 1, and this value is larger than 0.5. This phenomenon is ubiquitous for all type-II puzzles. We guess that CVX will transform some of the correct values 1 into the same value during optimization.

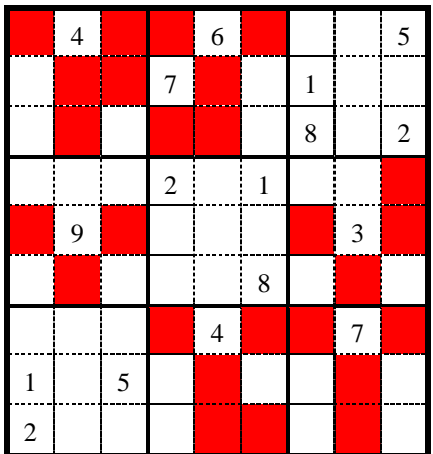

Fig. 3. Position of the key cells for a 17-clue puzzle

We also want to discuss our results in theoretical view. First of all, we think the key cell is an important and natural property for Sudoku puzzles. Then, what structure or property changes for filling the key cell? Can we find some rule to locate one key cell? The answer of the first question is related to the theoretical analysis for equivalence. We believe this is an interesting and useful open question for the theoretical study in CS. The study for the second question will produce some more effective algorithms to solve Sudoku puzzles.

## IV. IMPROVED ALGORITHM

In order to show our findings about the key cell is useful and natural, we propose an improved algorithm that uses an adopt threshold utilizing the phenomenon explained in the preceding section. The adopt threshold is set to the mode of the CVX solution. The algorithm framework can be described as follows:

**Improved Algorithm**:

(1) Solve the Sudoku puzzle using CVX

**if** "the result is not correct," **Do Once**

(2) Round the CVX solution $x$ to four significant digits. Extract the numbers from $x$ and obtain a set $S = \{x_i \,|\, 0.5 \leq x_i < 1\}$. Set the threshold value $t = \text{mode}(S)$

(3) Set the number to 1 in the CVX solution, whose value is equal to the threshold $t$. Fill the empty cells of Sudoku according these new 1s

(4) Solve the new Sudoku puzzle using CVX

The improved algorithm can solve 5923 (79.73%) type-II puzzles exactly, and the total accuracy rate is 96.94% for all the 17-clue puzzles. This result is a little better than the accuracy rate of the weighted L1 minimization algorithms (93%~94%) [21].

## V. CONCLUSIONS

In this study, the equivalence of L1 and L0 minimization is studied for a total of 17-clue puzzles. These puzzles are divided into two different types according to the solution uniqueness of their L1 minimization. For the type-I puzzles, the solution of L1 minimization is unique, and the sparse optimization algorithms can solve all of them exactly. By contrast, the solution of L1 minimization is not unique for the type-II puzzles, which cannot be solved by the algorithms correctly. Some cells of the type-II puzzles can change the equivalence of L0 and L1 minimizations. These phenomena may be helpful for the study of equivalence of L0 and L1 norm minimization in CS. We also propose an improved algorithm that uses an adaptive threshold according to the true value of the key cells. We may believe that the L1 minimization algorithms could be improved by suitable strategies through utilizing the prior for this problem and other L1 relaxation problems.

ACKNOWLEDGEMENTS

This work was supported by the National Natural Science Foundation of China (no. 61601518 and 61372172).